# Developing an Augmented Reality-Based Game as a Supplementary tool for SHS-STEM Precalculus to Avoid Math Anxiety

Carlo H. Godoy Jr.
*Technological University of the Philippines-Manila*
Manila, Philippines
carlo.godoyjr@tup.edu.ph

*Abstract*— **Math Anxiety is experienced by students. This is caused mainly by poor academic performance specifically in Calculus and Precalculus in Mathematics. From 2014 to 2016, an average rating of 70.33 percent for secondary education and 71.1 percent for elementary was generated during several exams, such as the National Achievement Test. The result produced is said to be below the national passing percentage, which shows that Filipino student performance in mathematics is very poor indeed. Taking into consideration the proof of the existence of these problems, it is necessary to try and test a need for a technology which has never before been used as a remedy. Augmented reality is deemed to be one of the advances in technology that is found to have a good effect not only on the health of the student, but also on the learning of Mathematics. Combining the two positive effects of Augmented Reality and the very impressive outcome of a game-based learning technique in mathematics, it is wise to say that an Augmented Reality Game-based learning method can be used as a supplement in Precalculus and Calculus teaching.**

## I. Introduction

An Augmented Reality or simply called AR is a system that allows an electronically-presented material to be introduced within the interaction of the user and the real world (Carbonell Carrera and Bermejo Asensio 2017). During the start of the new millennium, people are dreaming of this technology to be available since it was overwhelming on how this technology works during that era and no existing technology is comparable to this aside from fiction movies. As those movie introduces the nation with the technology, it became evident that people want this to be a reality and it happened (Peddie 2017). From that time, the desire of humanity for that technology gave birth to a new phenomenon that nobody knows will be popular for most of the people in the world. Nowadays, the technology emerged into different platforms and now it is being used as a supplement in different industries like healthcare, government, entertainment and education (Alkhattabi, 2017; Mahmood et al., 2017; Raiff, Fortugno, Scherlis, & Rapoza, 2018; Seevinck, 2017).

According to BusinessWorld (2019), it is very evident that there is a decline in the academic performance of students in the Philippines. Mobile game addiction of most Filipinos also post a great threat in the future of the society just like the declining rate of academic performance of the students (Fabito and Yabut 2018). In a study conducted by the National University of the Philippines, it is identified that mobile game addiction is a concern in the Philippines that grows rapidly which normally increases the risk in the physical as well as the psychological health of a certain individual (Fabito et al. 2018). Based on a report from UNESCO last 2008, the Philippines's Basic Education has low participation and achievement rates because it has fallen dramatically (Wilson Macha, Christopher Mackie 2018). Mobile game addiction and poor academic performance are among two of the most common problems a student in the Philippines is facing since mobile phones are known to be very affordable nowadays that leads the students to spend more time playing games or using social media than studying, especially if the subject is as hard as Precalculus and Calculus.

It is very evident that these problems are happening in the Philippines. There are students experiencing Math Anxiety. This is mainly caused by poor academic performance in Mathematics specifically in Calculus and Precalculus (ESTONANTO n.d.). During several examinations like the National Achievement Test, an average rating of 70.33% for secondary education and 71.1% for elementary was generated from the year 2014 to 2016 (Nepaya 2019). The generated result is said to be below the national passing percentage, which shows that the performance of Filipino students in mathematics is indeed very poor. On the other hand, mobile game addiction is very evident, as most of the places here in the Philippines are full of people playing Mobile Legends which is the most popular mobile game nowadays. In the past years, before the addiction in mobile legends arise, people walking in all areas of the Philippines always has phone cameras on as if something has been lost and needed to be found. It was then ruled out that what has been looking for is not visible with the naked eye but only with the cameras on due to an Augmented Reality game. This game was popularly known before as Pokémon Go. Due to this rapid growth of popularity, many gamers are suffering with different sickness as shown in all social media platforms. Many family members are posting these negative effects of the game to the gamer's social media accounts.

Poor Academic performance in college subjects like Differential Calculus, Integral Calculus and Solid mensuration on the other hand is found to be caused of weak background in Precalculus during the Basic Education years of the students (Agatep 2018). This poor academic performance with those mentioned subjects will lead to lower grades, social problems with performers in the class and in a worst-case scenario dropping out of the program. Both this Game Addiction and Poor Academic Performance in Mathematics, especially in Precalculus and Calculus has been a long problem of the people here in the Philippines. Several remedies in the past have been tried but until now these problems are still affecting the lives of every students here in the Philippines.

Considering the evidence of the existence of these problems, a need for a technology that has never been used before as a solution must be tried and tested. Augmented Reality or simply AR is considered as one of the advancements in the field of technology which is found to have a good effect not only with student's health but also provides interest in learning Precalculus ( Koivisto, 2019; Salinas, Patricia; Pulido, 2017). Combining the two positive effects of Augmented Reality and with the very promising result of a Game-Based learning method in mathematics, it is best to say that an Augmented

Reality Game-Based Learning approach can be made as a supplement in teaching Precalculus and Calculus (Tokac, Novak, and Thompson 2019). As known by many Filipinos, these two subjects require higher order thinking skills, which means that the students must have a healthy cognitive and emotional performance. According to Ruiz-Ariza et al, (2018) Pokémon Go has a very good effect on the cognitive function as well as the emotional well-being of adolescents. With, it is safe to say that an Augmented Reality Game, based on the contents of Precalculus and Calculus will be able to provide the same result. Thus, giving an advantage for the Philippines to have a healthy student with a high academic performance in Precalculus and Calculus.

## II. Related Literature and Studies

**Definition of Augmented Reality.** According to Rouse (as cited in Jung & tom Dieck, 2018) Augmented Reality or simply AR is the integration of information in digital format which includes live video on the real time environment of a certain user. In an augmentation of live videos, integrating a video picture to digital environment involves identification of an object replicated from the physical world features and will be captured as any format that will be considered as a video picture which will mean that increasing the responsiveness of the generated video picture to the state needed to control the object from the physical world itself (Kochi, Harding, Campbell, Ranyard &Hocking, 2017). In an augmented reality system, the integrated digital information can only be seen using a medium like phone cameras, but it will not be seen in the real world. These digital information can be represented in different forms like a stack of virtual cubes or manipulating a non-real object in many ways possible (Hilliges, Kim, Izadi, & Weiss, 2018).

**Effects of Augmented Reality in One's Health.** In line with the several Augmented Reality applications that has been created in the field of medicine, many positive effects have been studied and found out to be helpful in one's health. One of them is the positive effects provided by an Augmented Reality based application that enhances indoor navigation of the people with disability that uses a wheel chair (de Oliveira et al. 2016). Since individuals with reduced limb-related mobility impairments often face huge difficulties in participating in routine operations and moving around different settings, an application has been created to help them with their navigation. The novel characteristics introduced in the suggested architecture of the application, with unique emphasis on the usage of a Mobile Augmented Reality and its contribution as well as the capacity to define the best paths free of potential risks for wheelchair users, were able to provide important advantages for the indoor navigation of people with disabilities that uses a wheelchair compared to the current traditional techniques being used by the caregivers.

**Effects of Augmented Reality in One's Socialization.** Aside from the positive effects in the physical well-being of the players of pokemon Go, it also has a positive effect in the socialization skills of the players. Pokemon Go not only provides a novel gaming experience, but also can change how players perceive their physical realities. In relation to popular game-derived experiences and gratifications, (location-based) Augmented Reality Games can afford, for instance, outdoor adventures, community events, and health advantages, it provides a lot of exercise in socializing with other people (Hamari et al. 2018). Scientists and technicians rave about the capacity of various respondents to collaborate on a common project from distant places, educational experts are eagerly exploring methods to develop decentralized virtual classrooms, and social networking businesses such as Facebook state irrefutably that they see augmented and virtual reality as one of the main instruments that facilitate social interaction in the community (Metz 2016). Gamers are no different from these social media users. From forming an alliance to achieving a common goal to extremely competitive skill games, accommodating various users are just some of the significant factors in the growth of gaming software for augmented reality applications (Aukstakalnis 2017). In this manner, these gamers are not just forming an alliance in the game but also a blossoming that can extended not just inside the game but also in the physical world.

**Definition of *Game-Based Learning*.** The general definition given by education professionals as well as gamers that are also keen with their study habit, Game-based learning is an educational game that integrates learning content like theoretical and practical application of a subject to the idea of digital games to assist leaners, gamers and players learn about topics in relation to the subject matter. Prensky further described game-based teaching as a sort of game intended to balance the topic with gameplay and the player's capacity to maintain and apply the topic to the actual globe (Chang et al. 2017). Game-based learning incorporates tightly the learning content of a certain subject matter and the features of the computer games that enhance learning interests and motivations of the learners. In terms of leadership skills, Game-based learning (GBL) has gained attention from scientists and professionals as a universal strategy to building management abilities (Sousa and Rocha 2019).

***Augmented Reality for Mathematics Education***. An integrated STEM (Science, Technology, Engineering and Mathematics) lesson requires to participate and nurture students ' interest in real-world circumstances. While real-world STEM situations are naturally incorporated, the embedded STEM contents are rarely taught by school educators (Hsu, Lin, and Yang 2017). One of the hardest subjects of that track is Mathematics. One example of a Mathematics subject is Solid Geometry. To give a better experience in learning solid geometry, a study has been conducted to combine Augmented Reality (AR) technology into teaching operations designing a learning scheme that helps junior high school learners learn sound geometry (Liu et al. 2019; TeKolste and Liu 2018). Based on the result of the study, AR really gives a big leap in learning solid geometry. Another study deals with the use of AR in teaching and learning math that uses this technology to its complete benefit in providing concrete experience in interacting with revolutionary solids. At the end of the study, it was found out that Augmented Reality is beneficial in the understanding of computing solids of revolution volumes (Salinas and González-Mendívil 2017).

**III. PROPOSED METHODOLOGY**

    A. **Project Design**

Scavenger-Calc is a mobile game that uses Augmented Reality to tweak the images in the real life to have a virtual representation that can only be seen when you point your camera to a specific location. Scavenger-Calc will be designed similar to Pokemon-Go but with a twist of the traditional scavenger hunt game. This mobile game will be an investigatory game where in the learners will learn the module from beginning to end while solving the mystery of the game. If the mobile game will be designed in a somehow almost perfect manner, a game based learning approach for teaching Precalculus will control and make use of the student's intrinsic motivation together with their interest in playing and lead them in a very interactive way of solving exercises and seat works for the Precalculus subject.

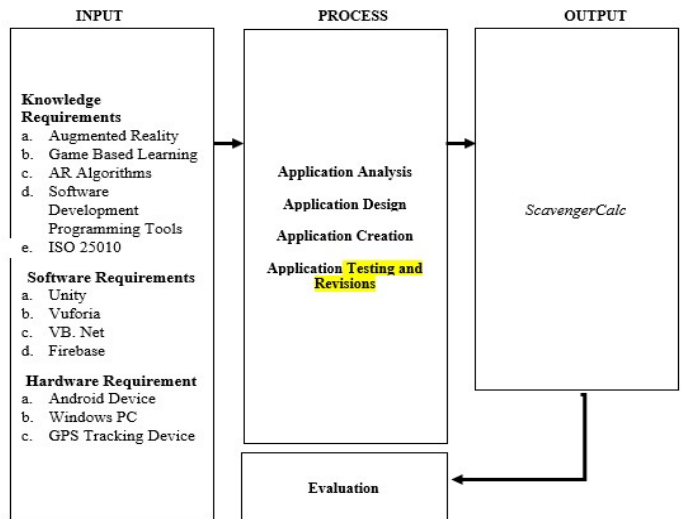

*Figure 1.* The Conceptual Model of the Study

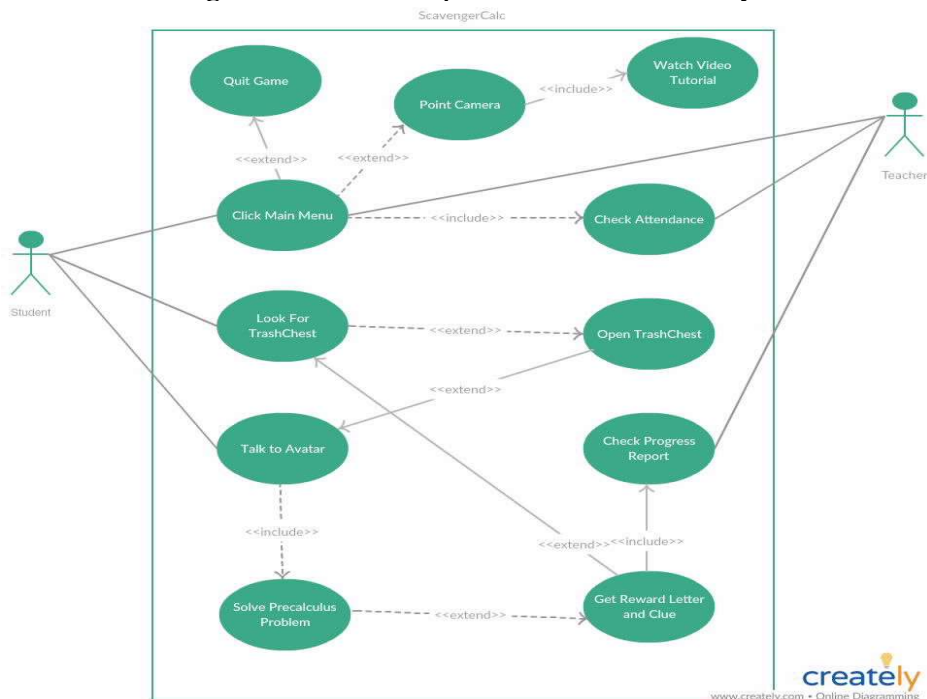

*Figure 2.* Use Case Diagram of ScavengerCalc

Figure 2 shows the Use Case Diagram of the application which is named as ScavengerCalc. This also illustrates and describes the interaction of the users with the system. The system has two types of users namely Student and the Teacher. Both the teacher and the student will have the ability to click on the main menu. Once on the main menu, the teacher can go ahead and check the attendance and the progress report of the students. On the other hand, the student will have two options after clicking on the main menu. The first option of the user is the tutorial mode. On the tutorial mode, the students can point out the camera in the physical environment.

Once the camera has detected any symbol or object related to any topic in Precalculus, it will go ahead and play the tutorial video for that certain topic. The second option for the student is to play the ScavengerCalc Mystery Game. The student will walk around the area and look for a TrashChest. Once the student found a TrashChest the option will be to open it or not. If the student chooses to not to open it, the student will just look for another TrashChest. If the studentchoose to open it, a Precalculus problem will appear and the student will need to solve the problem. If the student got the answer wrong, the student would need to walk around again looking for another TrashChest. If the student got the answer right, the student would get another the Reward letter and will have the chance to go on a quest to complete the Mystery word following the clue given by the Avatar. All the points garnered will be saved to the progress report afterwards.

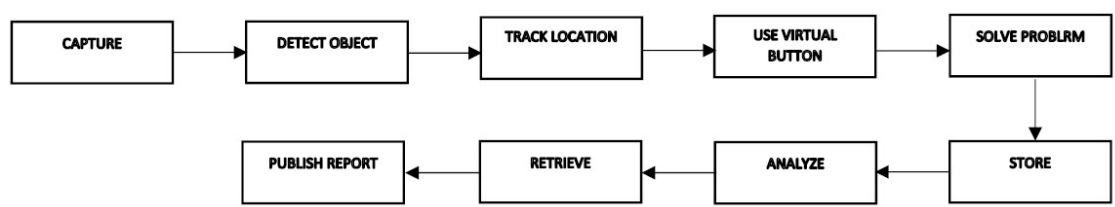

***Figure 3.*** Architecture Design of Location-Based Tracking Algorithm

Figure 3 illustrates the Architecture Design of Location-Based Tracking Algorithm. The process will start when the user pointed out the camera on the specified location while walking then when the camera pointed out the location where the TrashChest is located, it will prompt the user to open it or not. If the user decides to open the TrashChest and solve the problem with a correct answer, the score will then be generated as a data to be stored in the Cloud Database. When either the student or the teacher will check the generated result, the data analysis tool will then retrieve the data from the cloud database and generate a report from the retrieved data.

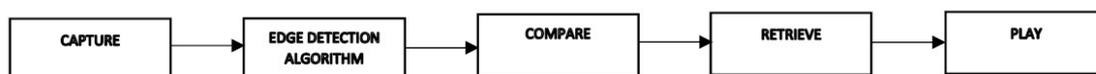

***Figure 4.*** Architecture Design of Image Recognition Algorithm

Figure 4 shows the Architecture Design of Image Recognition Algorithm. The process will start when the user points out the camera on the specified image. The image recognition algorithm will then read the image. The image will be check if it's a circle, parabola, hyperbola or elipse. Once recognized, image will then be evaluated. After evaluation, it will be compared to the videos stored in the multimedia database. Once the video related to the image has been identified, the video will then be retrieve and it will start playing as a tutorial video.

### B. Proposed Project Development Model

This research uses the Life Cycle of Software Development (SDLC). In each defined stages of the SDLC which will define the software development process, the framework will be described in every part of the process. It comprises of a comprehensive plan describing the development of the scheme.

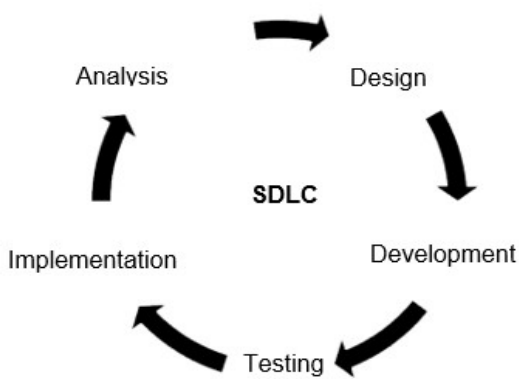

*Figure 5.* Software Development Life Cycle (SDLC) Model

The SDLC comprises of several clearly defined stages, namely Analysis, Design, Development, Test and Implementation..

A. Analysis
 1. Preparation of questions to be asked to the end-users and technical experts.
 2. Discussion with end users of like students, teachers and administrators.
 3. Discussion with technical professionals, comprising of developers of Android and IOS apps, to learn some best practices.
 4. Reviewing of existing application that has the same goal, issues, methodology and identify its weak points
 5. Definition of the process flow of the system using a diagram.
 6. Identification of the needed tools and software which will be needed for the development the application
B. Design
 1. List the system requirements and characteristics
 2. Build the following: Context Diagram, Data Flow Diagram (DFD), Students Low-Level Flowchart, Teachers Low-Level Flowchart, Use Case Diagram and Interface Design.
 3. Produce interface design layout and mock-up apps to list possible changes.
C. Development
 1. Use Vuforia and Unity to convert the final design layout to a working interface.
 2. Use MySQL to create the database and store processes.

      3. Use C # to create the API (Application Programming Interface) for the key system functionality.
      4. Use Vuforia and Unity to integrate the interface and API to create one working application.
  D. Testing
      1. Identification of tiny components of apps that are testable, called units.
      2. Test for correct operation of each unit.
      3. Test the system's general integrity.
      4. Test the system based on functionality, reliability, and compatibility using test instances and live testing with actual end-users.
  E. Implementation
      1. Request the application to be in Playstore and Apple Appstore.

### C. Operation and Testing Procedure

To ensure system quality, series of tests will be conducted for each module and by installing the application in a live environment, the system will be subjected to real-life testing.

**Functionality Test**. The functionality test will be used to guarantee that the application meets all the criteria and implements all the functionalities indicated in its functional requirements.

**Reliability Test**. This test will be performed to confirm that for a particular period the system can perform the designated tasks correctly. Live testing will be performed to check the system's reliability. The system will be implemented and exposed to real-world testing in a live setting.

**Portability Test**. Portability test will be used to explore the efficiency and effectiveness of the system with distinct devices that has distinct requirements such as variants of Android OS and IOS, hardware settings, screen sizes, and kinds or speeds of the network. The experiments will be provided in Table 8 to install and use the application in multiple kinds of setting. If all the anticipated behaviors are met, the testing is effective. The compatibility tests are provided in the Appendix for each module.

## IV. CONCLUSION AND FUTURE WORKS

Educational game serves as one of the applications of game-based learning method. Educational games were described as apps using video related mobile game features in producing interactive and immersive experiences in learning with specific objectives. These games provide challenge assignments, foster various interaction levels, and provide pleasant multimedia and immediate feedback (Denden, Essalmi, and Tlili 2017). Researchers have shown the efficacy of educational games in many fields such as language learning and mathematics. For future studies, an Augmented Reality Based mobile application named ScavengerCalc will be developed by the researcher as a supplementary learning tool in teaching and learning Precalculus. The proposed conceptual framework and the sdlc methodology discuss in this study will be used in the development to ensure the quality and correctness of the process. Once the mobile application has been developed, it must be evaluated by IT experts using ISO 25010. Based with their future